# An Almost Exact Linear Complexity Algorithm of the Shortest Transformation of Chain-Cycle Graphs

K. Yu. Gorbunov[1] and V. A. Lyubetsky[1]

**Abstract**
A "genome structure" is a labeled directed graph with vertices of degree 1 or 2. A set of operations over such graphs is fixed, and each of the operations has a certain cost, a strictly positive number. The transformation problem consists in the following: for given structures $a$ and $b$ and given costs, find a minimum total cost sequence of operations transforming $a$ into $b$ ("the shortest transformation of $a$ into $b$"). Each operation corresponds to an "event", the latter being a change in the graph caused by executing one of the operations over it. The possibility of assigning different costs is important in applications, since it allows to distinguish between frequent and rare events. Apparently, arbitrary costs make the problem NP-hard, which results in nontriviality of passing from one restriction on costs to another, if the problem is solved by a linear or at least polynomial algorithm (assuming that P≠NP). We propose a novel linear time and space algorithm which constructs a sequence of operations transforming $a$ to $b$ with total cost close or equal to the absolute minimum. Namely, if all the so-called DCJ operations have the same cost $w$ and if deletions and insertions have costs either both larger or both less than $w$, then the algorithm outputs a transformation of $a$ into $b$ with the total cost differing from the absolute minimum by at most $2w$ (in the former case) or equal to it (in the latter case). In some cases, the algorithm outputs an exact solution, e.g., in the case of circular genome structures. The condition on the costs of deletions and insertions can be omitted (although the proof described below does not include this general case).

**Keywords:** combinatorial optimization, exact graph algorithm, linear complexity

## 1 Introduction

Here and in what follows, a "genome structure" (we will simply say a *structure*) is any directed graph whose components are chains and cycles (including loops), and each edge has a name (which is a natural number). In other words, the degree of any vertex in a graph is either 1 or 2. The main theorem of the present paper says about structures with *all names different* in each individual structure; names in different structures may repeat. A structure has no vertices of degree 0, which will be referred to as *isolated vertices*. A usual set of operations over such graphs is fixed (see the next section), and each of the operations has a certain cost, a strictly positive number.

The problem of finding the minimum total cost sequence of operations transforming of one structure $a$ into another $b$ ("the shortest transformation") has been intensively studied in bioinformatics for a long time. Its occurrence is due to the fact that, at large evolutionary distances (globally, as they say), a genome in biology can be described as a labeled graph composed of chains and cycles (respectively, linear and circular chromosomes). For arbitrary operation costs, the problem is highly nontrivial mathematically because of NP-hardness, and thus is far from being solved. However, its importance in evolution research is widely recognized.

The problem traditionally uses some biological terminology. Namely, an edge with its name is called a *gene*, and a component of a graph, a *chromosome*. The graph (structure) itself is sometimes referred to as a *genome*. Such a global description of a biological genome, of course, ignores most of its real features; nevertheless, it happens to be very useful in studying many evolution properties of biological genomes. At the same time, a simplified description of a biological genome as a structure makes it possible to find efficient algorithms working with hundreds of genes simultaneously. Over the last 30 years, hundreds of papers on this problem were published, where it was usually considered in a biological context, whose overview can be found, e.g., in; Warnow [1]; Chauve at al. [2]; Pevzner [3]. We discuss this problem as pure mathematical, but in the next two paragraphs we briefly mention some of the most important original results in a biological context.

---

[1] Institute for Information Transmission Problems of the Russian Academy of Sciences,
E-mail: gorbunov@iitp.ru, lyubetsk@iitp.ru

In Sankoff et al. [4], genomic rearrangements in ancestors of twenty mitochondrial genomes were analyzed. In Hannenhalli and Pevzner [5] there were considered genomes with equal gene content (i.e., the sets of names in genomes *a* and *b* coincide), a single chromosome, and reversal as an operation over the genome. The runtime of the obtained algorithm is estimated by a fourth-degree polynomial. The *distance* between two genomes is understood as the number of applications of the given operations in the shortest transformation between them. In Hannenhalli and Pevzner [6] there was considered the problem of finding the distance between the same genomes as in Hannenhalli and Pevzner [5] but with a larger list of operations: fusion, fission, reversal, and translocation. The runtime of the obtained algorithm was also a fourth-degree polynomial. The above operations are expressed through the so-called DCJ operations, whose definitions are recalled in the next section. In those papers, the problem was considered under the equal gene content and without costs of operations, and simply the number of operations in a transformation was minimized. In Alekseyev and Pevzner [7] (for equal gene content) there was proposed an important construction of a breakpoint graph, and the transformation problem was reduced to the problem of transforming the breakpoint graph to a final form (these definitions see below).

In Braga et al. [8]; Compeau [9], Silva et al. [10], Compeau [11], Silva et al. [12] this problem was studied under unequal gene content and unequal operation costs. These papers describe linear time algorithms constructing the shortest transformation of one structure into another, assuming equal costs of DCJ operations and, independently, equal costs of insertions and deletions (referred to as *indels*) of a gene. In the present paper we also assume equal cost $w$ of DCJ operations, but costs $w_d$ and $w_i$ of deletion and insertion can be different, though both greater or less than $w$. Algorithms and proofs from the papers mentioned above cannot be extended to these cases. Note that possibility to take different costs $w_d$ and $w_i$ is essential in some applied and in particular biological contexts. For instance, mitochondria of all vertebrates have the same set of genes and each gene has a unique function. However, the order of these genes can vary among animals: it is substantially different in the purple sea urchin *Strongylocentrotus purpuratus* than in vertebrates (Jacobs et al. [13]). Thus, gene losses and acquisitions are much rare here than DCJ operations, and rare events should be assigned higher costs. A similar pattern is observed in plastids of certain plants; e.g., the order of plastid genes in the red alga *Porphyridium purpureum* notably differs from that in other red algae. In the next case where one genome is smaller than the other, the loss cost should exceed that of acquisition. For example, the genome size amplifies after hybridization of two species as in the case of polyploid strawberry or wheat (Bors and Sullivan [14]). Polyploidy is much more common in plants than in animals. Among dioecious animals, it was described in nematodes including ascarids and certain amphibians. However, other full-genome duplications have been recently reported in chordates (Putnam et al. [15]). The genome reduction can also occur; e.g., many genes (paralogs) are lost after a full-genome duplication.

Algorithmic and mathematical hardness of passing from one condition on the costs to another follows from the fact that the problem with arbitrary costs is apparently NP-hard and therefore definitely cannot be solved even by an algorithm of polynomial complexity (assuming that P≠NP). In the present paper we consider the problem under *unequal gene content* and minimize by a linear algorithm the sum of costs of all operations applied in a transformation. This sum is referred to as the *total cost* of a transformation.

Let us recall related problems, though they are not considered here. If edge names in a structure may repeat (paralogous genes), then this problem of the shortest *transformation with paralogs* was investigated in Shao et al. [16] under the following constraints: equal gene content, equal number of the same name paralogs in structures, and DCJ operations only. In Martinez et al. [17] the problem was analyzed with DCJ operations and unequal gene content, but paralogs that do not match any paralog in another structure are discarded. In Lyubetsky et al. [18]; Lyubetsky et al. [19] the general case of the transformation with paralogs problem was reduced by a linear algorithm to a special case of an integer linear programming problem. In Lyubetsky et al. [18–19] there was also considered the general case of the problem of the shortest *reconstruction* of structures (with or without paralogs) defined at leaves of a tree, i.e., the problem of extending structures from leaves to inner nodes of the tree. The tree itself can be either given or found in the course of the solution.

# 2 Definitions and preliminaries
## 2.1 Equal gene content and equal operation costs
In algorithmic considerations, names of genes are natural numbers. *Equal gene content* means that the sets of names in *a* and *b* coincide.

DCJ operations over a structure are well known, see Bergeron et al. [20]. These are the following. *Double intermerging*: cut two vertices and merge the four ends thus formed in a new way; Fig. 1a. S*esquialteral intermerging*: cut a vertex and merge one of the newly formed ends with any free end; Fig. 1b. *Single merging* and *cut*: merge two free ends or cut a vertex; Fig. 1c.

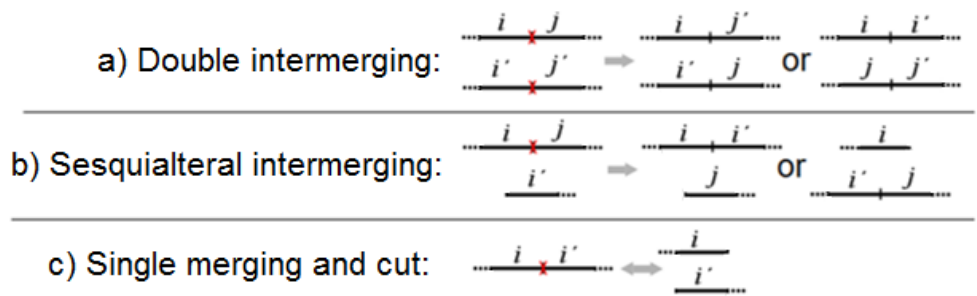

**Figure 1.** DCJ operations over a structure (edge directions are not shown).

Under equal gene content, the *transformation problem* consists in the following: for given structures *a* and *b*, find a minimum-length sequence of operations transforming *a* to *b*. For solving this problem, of major importance is the notion of a *breakpoint graph a+b* for the structures *a* and *b*, which was proposed in Alekseyev and Pevzner [7]. This is an undirected loopless graph whose vertices are extremities of all genes in *a* and *b*. Let us denote the beginning of a gene with name $k$ by $k_1$, and its end, by $k_2$. By definition, an edge connects two vertices in *a+b* if the corresponding extremities are merged in *a or* in *b*; the edge is labeled by *a* or *b*, respectively. Thus, *a+b* contains information on merged extremities in *a* and *b* simultaneously. Components of *a+b* are chains and cycles with alternating *a*- and *b*-edges. Isolated vertices may occur; they correspond to extremities that are not merged with any others.

The graph *c+c* consists of isolated vertices and cycles of length 2; it is said to be *final* (or: in *final form*). It is easily seen that the transformation problem is equivalent to the problem of reducing *a+b* to a final form by operations which are analogs of operations over pairs of structures, see Alekseyev and Pevzner [7]. Namely, an operation is performed over one element of a pair of structures only, and the second is unchanged; after that, an operation over *a+b* is defined by commutativity of the corresponding diagram for any pair <*a,b*>. These operations are as follows. *Double intermerging* (Fig. 2a): deleting two equally labeled edges and connecting the four ends thus formed by two new edges with the same label. *Sesquialteral intermerging* (Fig.

2b): deleting an edge and joining one of its ends with a free vertex that is not incident to an edge with the same label by an edge with the same label. *Single merging* (Fig. 2c, right to left): adding an edge (say with label *a*) between two vertices such that each of them is not incident to an edge with label *a*. *Cut* (Fig. 2c, left to right): deleting any edge. As a result, labels of edges in *a+b* alternate.

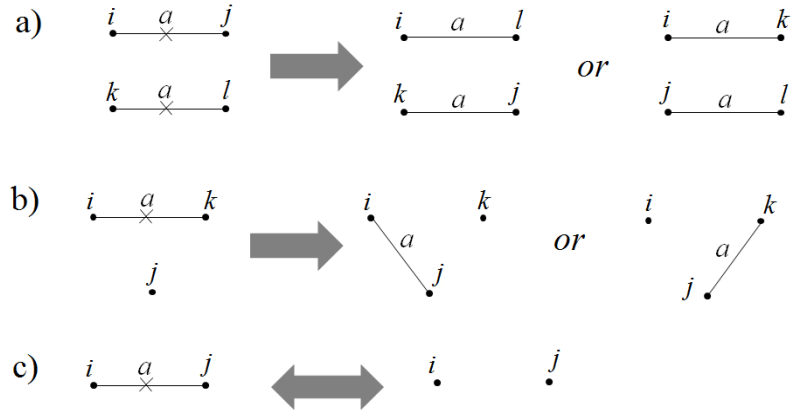

**Figure 2.** Operations over the breakpoint graph *a+b* in the case of the equal gene content of structures *a* and *b*. Similarly for a *b*-edge.

**2.2 Unequal gene content and unequal operation costs**
Recall that we always seek for a transformation of a structure *a* to a structure *b*. *Throughout what follows*, we consider *unequal gene content*, i.e. sets of names in *a* and *b* need not coincide. A gene that occurs in both *a* and *b* is said to be *common*; a gene occurring in only one of the structures is said to be *special* (*a*- and *b*-special).

In this case, besides DCJ operations, which will be referred to as *standard*, we need *supplementary* operations over a structure: *deletion* and *insertion* of genes. Namely, we can delete a connected fragment of *a*-special genes and insert a connected fragment of *b*-special genes; if a deleted fragment has adjacent genes, merge their extremities. Such fragments will be referred to as *a*- or *b*-segments, respectively. An inclusion-maximal segment will be called an *a*- or *b*-block.

If each operation is assigned a cost, a strictly positive number, then the *transformation problem* is naturally generalized: find a transformation of *a* to *b* with the minimum total cost (say, *minimum cost*). Such a sequence will be referred to as *shortest*.

The result of the paper is the following

**Theorem 1.** *We describe a linear time and linear space algorithm for the transformation of any structures to each other. If standard operations have equal costs w and the deletion and insertion operations have costs both larger or less than w, then the algorithm outputs a transformation of a into b with total cost differing from the absolute minimum by at most 2w (in the former case) or equal it (in the latter case).*

In Gorbunov and Lyubetsky [21] it was proved that among shortest sequences there is one in which all deletions precede all insertions; therefore, the problem is equivalent to the following: transform *a* and *b* separately into any structure *c* by two sequences of operations (now without deletions) with the minimum total costs of these sequences. Thus, as in the preceding subsection, we look for a structure *c* to which *a* and *b* are independently transformed by shortest (in aggregate) transformations. To this end, we extend the notion of the graph *a+b*. Specifically in the case of unequal gene content, the breakpoint graph *a+b* has the following vertices: *conventional*, which are extremities of common genes, and *singular*, which correspond to *a*- and *b*-blocks and are labeled as *a*- or *b*-vertices. Each block is accompanied with a sequence of names for genes forming the block. A conventional edge in *a+b* connects conventional vertices as defined above. A singular edge in $a+b$ connects a conventional vertex with a singular one if an extremity in $a$ or $b$ corresponding to a conventional vertex is merged with an extremity of the corresponding block ($a$- or $b$-edge). A conventional isolated vertex is an extremity of a common gene that is merged in neither $a$ nor $b$. A singular isolated vertex is a linear block. A loop is drawn at a singular vertex that corresponds to a cycle. A *hanging* edge is a singular edge incident to a singular vertex of degree 1. The graph $a+b$ is undirected; its components are chains, isolated vertices, cycles, and loops. Non-hanging singular edges occur in $a+b$ in pairs — these are edges incident to one singular vertex; it is convenient to regard such a pair as a double edge, which allows us to preserve the alternating property of $a$- and $b$-edges. The *size* of a component in $a+b$ is the number of conventional edges plus half the number of singular non-hanging edges in it. For conventional isolated vertices and loops, the size is defined to be 0, and for singular isolated vertices, −1. The definition of a final graph (final form) remains the same, this is a graph of the form $c+c$. Fig. 3 shows a typical example of an $a+b$, in Fig. 4a we give an example of structures $a$ and $b$, and Fig. 4b shows their breakpoint graph $a+b$. A pair of structures $a$ and $b$ is transformed by a linear algorithm to the breakpoint graph $a+b$, and vice versa.

Other (though close in ideas) definitions of a graph combining the information about original structures $a$ and $b$ in the case of unequal gene content were proposed in Yancopoulos and Friedberg [22] and later in Braga et al. [8] and Gorbunov and Lyubetsky [21].

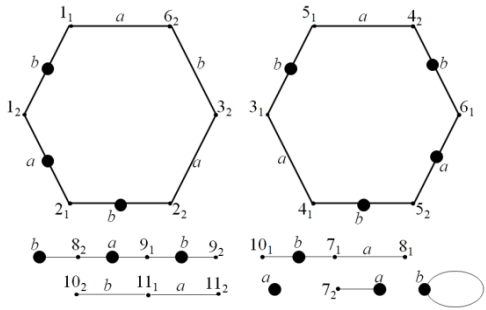

**Figure 3.** Example of a breakpoint graph $a+b$ for structures $a$ and $b$ with unequal gene content, which can be reconstructed given $a+b$. Conventional vertices are shown as small circles, and singular vertices, as big circles.

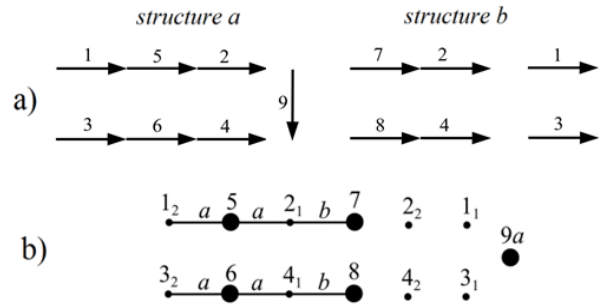

**Figure 4.** (a) Genome structures $a$ and $b$ with unequal gene content, and (b) their breakpoint graph $a+b$. Conventional vertices are shown as small circles, and singular vertices, as big circles.

We define five operations over $a+b$. These are as follows. *Double intermerging* (DM) is defined as above (Fig. 5a): if an edge with singular endpoints is formed, both $a$- or both $b$-, then the edge is replaced with one singular vertex, which is assigned the concatenation of sequences of the two original singular vertices («*joining*» of singular vertices); if an operation involves a loop, the loop vertex is considered as having two ends. *Sesquialteral intermerging* (SM, Fig. 5b): deleting an edge and drawing a new edge with the same label from one of newly formed free endpoints to a free conventional vertex which is not incident to an edge with the same label or to a singular vertex of a hanging edge or a singular isolated vertex (i.e., of degree ≤1) with the same label (along with possible *joining* of two singular vertices). *Single merging* (**OM**): adding an $a$-edge between free vertices each of which is a conventional vertex (joined to a $b$-edge) or an isolated vertex, an $a$-hanging vertex or an $a$-singular isolated vertex; along with possible joining of two singular vertices (Fig. 5c, right to left and Fig. 5c'). The same for $a$ and $b$ interchanged. *Cut* (Cut, Fig. 5c, left to right): deleting any edge. We define a supplementary operation: *removal* (Rem) of a singular vertex (Fig. 5d, d', d''). Specifically, if it is of degree 2, it is removed, and edges incident to it are merged into a single edge with the same label; if it is of degree 1, it is removed together with the edge incident to it; if it is of degree 0 or has a loop, the vertex and the loop are removed. Reduction of the graph *a+b* to a final form by definition uses the above-described operations.

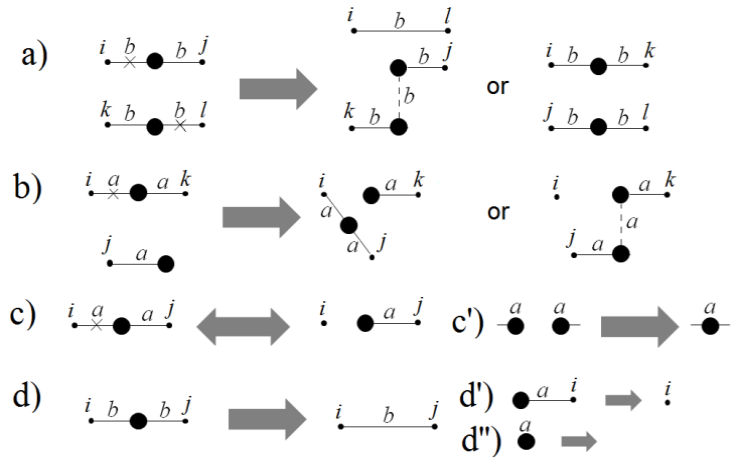

**Figure 5.** Operations over the breakpoint graph *a+b* in the case of unequal gene content of structures *a* and *b*.

*Costs of operations* over any breakpoint graph *G*, in particular *a+b*, are defined according to original costs. Specifically, the cost of removal of an *a*-vertex equals the original cost of deletion; that of removal of a *b*-vertex equals the original cost of insertion; of a cut of an *a*-edge, the original cost of a cut; of a cut of a *b*-edge, the original cost of a single merging; of a OM operation by an *a*-edge, the cost of a single merging; of a OM operation by a *b*-edge, the original cost of a cut; the costs of DM and SM equal the original costs of the respective operations. The *minimum cost* of *G* is defined to be the minimum total cost (under the above-defined costs) of a transformation (called *reduction*) of *G* to a final form. Transformation problem consists in finding a transformation itself and its minimum cost.

Recall that linear equivalence of two problems means the existence of a linear complexity algorithm transforming the problems to each other. Theorem 2 below literally repeats Corollary 5 from Gorbunov and Lyubetsky [21], whose proof does not use the equality of all operation costs.

**Theorem 2.** *Let the costs of all standard operations be equal. The problem of the shortest transformation of a to b is linearly equivalent to the problem of reduction the graph a+b to a final form. The minimum costs of both transformations are equal.* □

## 3 Description of the algorithm

We solve the problem of transforming *a* to *b* with arbitrary costs as the problem of reducing *a+b* to a final form with the costs specified above, which assumes that all *DCJ operations have equal costs* (the condition of Theorem 2). We solve this problem by a linear complexity algorithm whose almost exactness is proved in Section 4 below **provided that** $\min\{w_1,w_2\} \geq w$ (the condition of Theorem 3). The **second case of Theorem 1**, where $\max\{w_1,w_2\} \leq w$ (then the increment is 0, which implies that the algorithm attains the exact minimum) is proved similarly.

We may divide all costs by *w*; therefore, we assume that *w*=1, and it suffices to prove the exactness under the condition $1 \leq w_1 \leq w_2$, since the symmetric case is completely analogous. Note that we have also tested an implementation of the algorithm with randomly chosen costs, and it was found to show high performance; this separate subject is not presented in this paper. However, a computer program realizing the algorithm — together with biologically meaningful examples of its application — is available at http://lab6.iitp.ru/ru/chromoggl/ (in Russian).

The linearity of the algorithm will be evident, in contrast to its exactness. The main idea of the latter one consists in using an easily computable characteristic of any breakpoint graph *G*, a natural number *t(G)* (see Claim 2), which, roughly speaking, decreases by 1 at each step of the algorithm and equals 0 at a final graph; then the number of algorithm steps is *t(G)*. More precisely, for any graph $G$ and any operation $o$, we have the «*t*-property»: $t(G)-t(o(G)) \leq 1$ (see Claim 1). The other idea ensures that the absolute minimum of the total cost (called the *minimum cost*) of reducing an arbitrary breakpoint graph $G$ – in particular, $a+b$ – to a final form is close to the total cost of the sequence of operations output by the algorithm when applied to $G$. Specifically, we denote the first of them, the *minimum cost*, by $C(G)$, and the second, the *algorithm cost*, by $c(G)$. Then any operation $o$ applied to $G$ reduces $c(G)$ when passing from $G$ to $o(G)$, roughly speaking, by at most $c(o)$, the *cost of the operation $o$*. In other words, $c(G)-c(o(G)) \leq c(o)$. Thus, $t(G)$ and $c(G)$ play the role of a kind of a complexity estimate for $G$. We further formalize these ideas in Claims 1-2 and Lemma 2, respectively.

The algorithm consists of 8 stages, the 2$^{nd}$ and 3$^{rd}$ of them containing numerous substages. Anticipating things, by an *autonomous reduction* of any breakpoint graph $G$, in particular $a+b$, we call a sequence of operations consisting only of Stages 1, 4, 5.2, and 8 defined below. In essence, the autonomous reduction consists in cutting off all conventional edges (Stage 1), circularizing all chains (Stage 4), reducing cycle sizes (Stage 5.2); and circularizing a simplest chain of the form *aa* or *bb* and removing all singular vertices (Stage 8).

Now we start the description of our algorithm. At **Stage 1**, from all components except cycles of length 2, conventional edges are cut off, which are circularized into cycles of length 2. Namely, *cutting off* consists in a DM of a pair of adjacent edges or in SM of an adjacent edge, if the second adjacent edge is absent. At **Stage 2** sequences of operations (named *interactions*) are performed. The interactions are performed *over chains only* and reduce their number, as well as, in a sense, complexity of $G$. The interactions are performed in the order listed below and use the notion of a *type of chain*. At all stages (except for Stage 7), each operation transforms the graph $G$ into a graph $o(G)$ with a strictly lower value of $t$.

A chain of an odd (even) size in $G$ is said to be odd (respectively, even); zero is an even number. Hereinafter, the type of chain is *defined* after cutting off all conventional edges, and does not depend on the order of their removing (Gorbunov and Lyubetsky [21], Lemma 6). But the definition of a type is applicable to components of $G$ that contain arbitrarily many conventional edges, as well as vertices. *Type 1a* is an odd chain with one hanging *b*-edge; *2a* is either an odd chain with two hanging *b*-edges or a *b*-singular isolated vertex (types *2a\** and *2a'*); *3a* is either an odd chain without hanging edges and with two extremity *a*-edges, or a singular *a*-vertex connected with exactly two conventional isolated vertices (types *3a\** and *3a'*), Fig. 6. Respective *b*-types are defined similarly to *a*-types. Next, *type $1_a$* is an even chain with one hanging *a*-edge of size either strictly greater than 0 or exactly 0 (types $1_a*$ and $1_a'$); 2 is an even chain with two hanging edges of size either strictly greater than 0 or exactly 0 (types 2* and 2'); 3 is an even chain without hanging edges of size strictly greater than 0; 1 is the union of types $1_a$ and $1_b$; 0 is a chain without singular vertices.

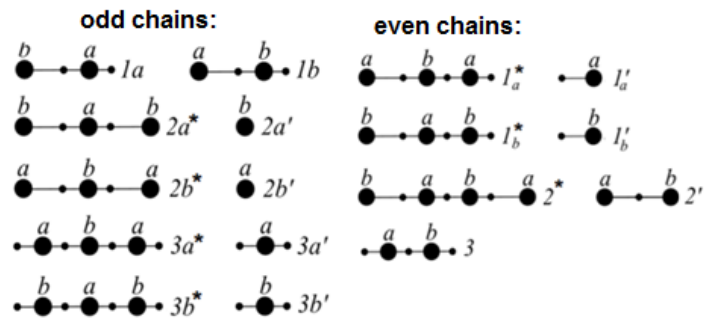

**Figure 6.** Possible types of chains in a graph $G$.

Components of $G$ are divided into *types* (*a,b*)-, *a*-, *b*-, and 0- according as they contain both *a*- and *b*-vertices, *a*- but no *b*-vertices, *b*- but no *a*-vertices, or no singular vertices at all. *Type 2c* is an (*a,b*)-cycle of size 2 with exactly one *a*-vertex and exactly one *b*-vertex.

7   --

Any *interaction* means performing the corresponding operations 1–3 applied to chains on the left-hand side of its equality and resulting in chains on the right-hand side of its equality. Interactions are performed according to the parentheses in them. Thus, interactions of **Stage 2** are the following: $1_b'+2a'=2a'$, $1_b'+3b'=3b'$, $1_b'+1_b'=1_b'$ (Fig. 7) and symmetric to them (SM operations); $1a+1b=1_b^*$, $2b+3a=1_a$, $2a+3b=1_b$, $2+3=1_b^*$ (SM operations), $(1a+2b)+3=1_b^*$, $(1b+2a)+3=1_b^*$, $(1b+3a)+2=1_b^*$, $(1a+3b)+2=1_b^*$ (SM), $1a+2=2a^*$, $1b+2=2b^*$, $1a+3=3a^*$, $1b+3=3b^*$ (SM); $(1a+2b)+(1a+3b)=1_b^*$, $(1b+2a)+(1b+3a)=1_b^*$ (SM), $(1a+1a)+2b=2a^*$, $(1b+1b)+2a=2b^*$, $(1a+1a)+3b=3a^*$, $(1b+1b)+3a=3b^*$ (OM and SM operations), $1a+1a=3a^*$, $1b+1b=3b^*$ (OM operations), $1a+2b=2$, $1b+2a=2$, $1b+3a=3$, $1a+3b=3$ (SM operations), $(2a+(2b+3))+3=1_b^*$, $((3a+(3b+2))+2=1_b^*$ (SM), $(3a+2)+2=2a^*$, $(3b+2)+2=2b^*$, $(2a+3)+3=3a^*$, $(2b+3)+3=3b^*$, $2a+(2b+3)=2^*$, $3a+(3b+2)=3$ (SM). «Symmetric» means interchanging *a* and *b*. Any interaction and any stage of the algorithm are performed repeatedly while it can be applied.

More precisely, stages of the algorithm assume choosing between several variants. For example, at Stage 1 this is choosing a sequence in which conventional edges are cut off. At Stage 2 interactions are described as operations over types rather than over breakpoint graphs *G* themselves. The output of the algorithm depends on these choices. However, as we will show, each of such results is a shortest reduction of the original graph *G* to a final form.

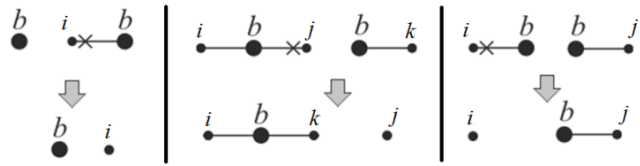

**Figure 7.** Interactions $1_b'+2a'=2a'$, $1_b'+3b'=3b'$, $1_b'+1_b'=1_b'$ (SM operations).

In (Gorbunov and Lyubetsky [21], Theorem 6), it is proved that application of Stages 1–2 to the graph *G* followed by any autonomous reduction results in a final graph with the *minimum number* of operations, which will be *denoted* by $t(G)$. This immediately implies the following.

**Claim 1**. *For any graph G and any operation o, we have «t-property»*: $t(G)–t(o(G)) \leq 1$. □

In fact, a stronger claim $|t(G)–t(o(G))| \leq 1$ holds, which we do not use here.

At **Stage 3** the following interactions are performed, for which we do not list symmetric interactions and resulting cycles of length 2, as well as conventional isolated vertices:
1) *a*-loop + «component with a singular *a*-vertex» = the same component (DM or SM operations);
2) $2'+2'=2a^*$ (OM operation);
3) $2a'+2a'=2a'$, $2a'+2a^*=2a^*$ (OM);
4) $3a'+3a'=3a'$, $3a'+3a^*=3a^*$ (OM operations and cuttings-off of a conventional edge);
5) $3a^*+3b'=3$, $3a'+3b'=3$ (SM operations); the proof of exactness becomes easier if the second interaction is performed only in the presence of $(a,b)$-components;
6) $2a'+2b^*=2+1_a'$ (SM; Fig. 8) and $2+1_a'=2$, which is defined in item 19 below;

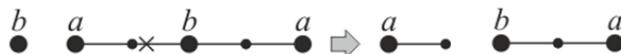

**Figure 8.** First operation (SM) of the 6th interaction: $2a'+2b^*=2+1_a'=2$.

7) $3a'+2'=1a$, $3a'+2^*=1a$ (SM);

8 --

8) $2a'+3=1a$ (SM);
9) $1a+3a'=1a$ (OM operations and cuttings off a conventional edge);
10) $1a+2a'=1a$ (OM);
11) $3a'+3=3$ (OM and cutting off a conventional edge);
12) $2a'+2'=2'$, $2a'+2^*=2^*$ (OM);
13) $1_a'+1_a^*=1_a^*$ (SM);
14) $1a+1_b'=1a$ (SM);
15) $1a+1_a'=1a$ (SM);
16) $2a'+1_b^*=2a^*$, $2a^*+1_b'=2a^*$ (SM);
17) $3a'+1_a^*=3a^*$, $3a^*+1_a'=3a^*$ (SM);
18) $3+1_a'=3$ (SM);
19) $2'+1_a'=2'$, $2^*+1_a'=2^*$ (SM).

Chains of types $2a'$, $2b'$, $3a'$, $3b'$, $1_a'$, and $1_b'$ are called *problem chains*; they have a single *a*- or a single *b*-vertex. Note that after Stage 3, there remain at most one problem chain with a single *a*-vertex and at most one with a single *b*-vertex.

Let $w_d=w_a$ be the cost of removal of an *a*-vertex, and $w_i=w_b$ be the cost of removal of a *b*-vertex.

At **Stage 4** chains of strictly positive sizes are circularized. Specifically, if a chain is odd and its ends are *a*-conventional or *b*-singular, then we perform an OM, joining these ends by a *b*-edge (and similarly for *a* and *b* interchanged). If it is even, we perform an SM, removing the external edge with a conventional end (if it exists) or the second edge from the end. In the first case, a conventional isolated vertex is formed and an edge with the same name as the deleted one is drawn to the opposite end. In the second case, the external edge forms a chain of length 1, and an edge with the same name as the deleted one is drawn. The proof of exactness becomes easier with the following specification: chains of type $3a'$ are circularized into cycles only if $w_a \geq 2$, and similarly for $w_b$.

At **Stage 5** the following interactions are performed:
1) $(a,b)$-cycle+$(a,b)$-cycle = $(a,b)$-cycle + cycle of length 2 (DM with joining of two singular vertices; Fig. 9);
2) $(a,b)$-cycle = $(a,b)$-cycle of a smaller size + cycle of length 2. This interaction consists of two operations: DM, which reverts a fragment of the cycle with joining of singular vertices, and cutting off the thus formed conventional edge; Fig. 10.

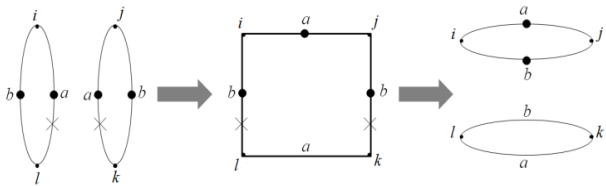

**Figure 9.** First interaction of Stage 5: merging two $(a,b)$-cycles into one $(a,b)$-cycle with a conventional edge followed by shortening of the $(a,b)$-cycle with joining of singular *b*-vertices (two DM operations).

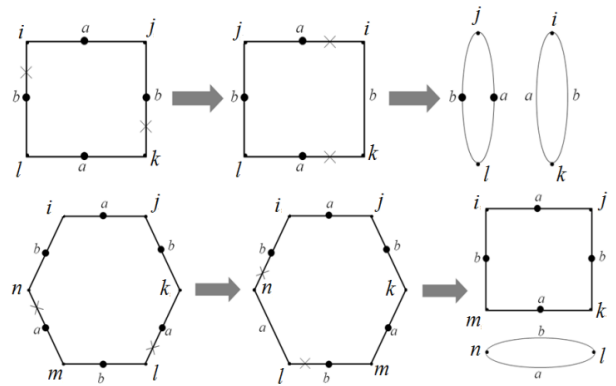

**Figure 10.** Second interaction of Stage 5: detaching a cycle of length 2 from an $(a,b)$-cycle of size 4 (top) and of size 6 (bottom). The first operation results in a conventional edge. The same for other sizes.

At **Stage 6** the following interactions are performed:
1) $2c+2a'+2b'=2'$ (two SM operations); Fig. 11.
2) $2c+2'=2'$ (DM and cutting off a conventional edge); Fig. 12. After Stage 6, at most one $(a,b)$-component remains.

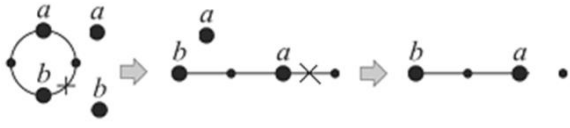 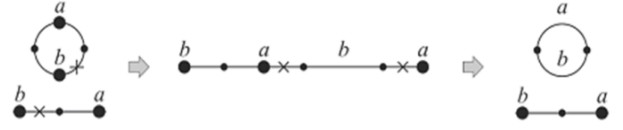

**Figure 11.** First interaction of Stage 6: $2c+2a'+2b'=2'$.　　**Figure 12.** Second interaction of Stage 6: $2c+2'=2'$.

The next stage, **Stage 7**, is applied only if $w_b>2$. At this stage $a$- and similar $b$-interactions are performed; moreover, the $a$-interactions are performed only if $w_a>2$.

1) $2c$-cycle+$a$-cycle = $2c$-cycle + cycle of length 2 (two DM operations); Fig. 13.
2) $a$-cycle+$a$-cycle = $a$-cycle + cycle of length 2 (two DM operations); Fig. 14.

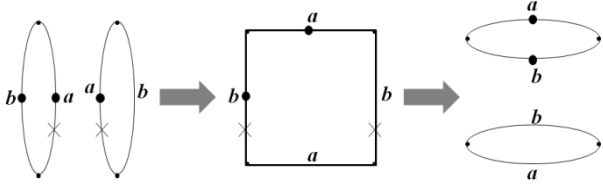 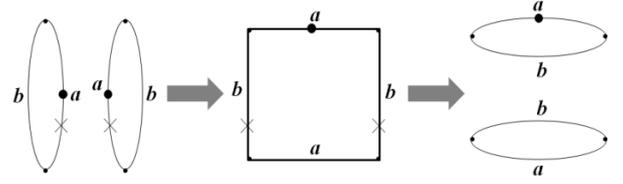

**Figure 13.** Stage 7, the 1st interaction: removal of an $a$-vertex is replaced with two intermergings, which is advantageous in weight. Similarly to the operations shown in Fig. 9, but $b$ in the premise is a conventional edge.

**Figure 14** Stage 7, the 2nd interaction: two $a$-cycles. Similarly to the 1st interaction, but with two conventional $b$-edges in the premise.

3) $a$-cycle+$2'$-chain = $2'$-chain + cycle of length 2 (SM and cutting off an external conventional edge); Fig. 15.
4) $2c$-cycle+$2b'$-chain = $2c$-cycle (SM and OM); Fig. 16.

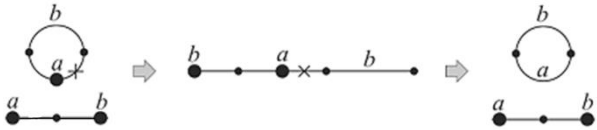 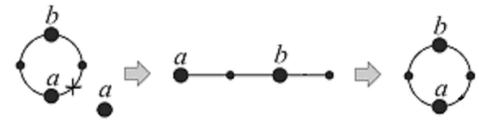

**Figure 15.** Stage 7, the 3rd interaction: $a$-cycle and chain $2'$.

**Figure 16.** Stage 7, the 4th interaction: $2c$-cycle and chain $2b'$.

5) $2c$-cycle+$1_a'$-chain = $2c$-cycle (two SM operations); Fig. 17.

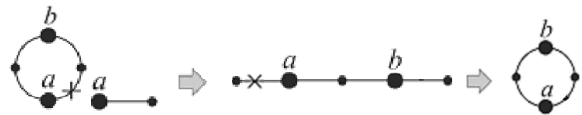

**Figure 17.** Stage 7, the 5th interaction: $2c$-cycle and chain $1_a'$.

6) $a$-cycle+$2b'$-chain = $2b'$-chain + cycle of length 2 (SM and cutting off an external conventional edge); Fig. 18.
7) $a$-cycle+$1_a'$-chain = $1_a'$-chain + cycle of length 2 (SM and cutting off an external conventional edge); Fig. 19. And symmetrically if $w_a$ and $w_b$ are interchanged.

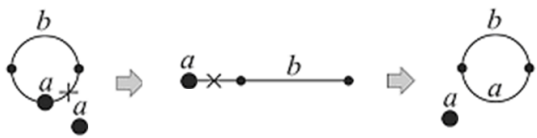 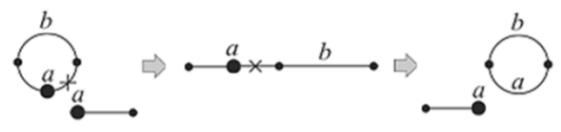

**Figure 18.** Stage 7, the 6th interaction: $a$-cycle and chain $2b'$.

**Figure 19.** Stage 7, the 7th interaction: $a$-cycle and chain $1a'$.

After Stage 7, if it is invoked, there remains at most one component containing an *a*- (if $w_a>2$) or *b*-vertex (if $w_b>2$).

At **Stage 8**, chains $3a'$ and $3b'$ are circularized, and singular vertices are removed.

## 4 Exactness of the algorithm

An operation applied to a breakpoint graph *G* is said to be *special* if performing it strictly reduces the number of singular vertices. Such operations are, first, the DM, SM, and OM operations, which include *joining* of singular vertices, and second, the Rem operation itself. Note that the DM and SM occur at Stages 1–7 only, and Rem, at Stage 8 only.

*Let S* be the number of nonspecial operations in Stage 1, and *D* the number of chains of types 1*a*, 1*b*, 3*a*, 3*b*, and 3. Here is an important notion: *P* is the difference between the *number* of operations in autonomous reduction of *G* and the *number* of operations in the transformation obtained after executing Stages 1–2 followed by autonomous reduction. Thus, the only difference is that the subtrahend also includes Stage 2. Clearly, $P(G_1)=0$ for any graph $G_1$ resulted from Stage 2.

Let $C_a$ and $C_b$ be the *numbers* of *a*- and *b*-cycles (not loops) in *G*; and let *B* be the *number* of singular vertices in *G*. An indicator is a function that equals 1 if a fixed property of a graph *G* is satisfied and equals 0 otherwise. Introduce the following *notations*: $I_{ab}$ is the indicator of an (*a*,*b*)-component in it; $I_{pa}$ the indicator of an *a*-loop if there are no other components with an *a*-vertex (the same for $I_{pb}$); $I_{ca}$ the indicator of an *a*-cycle (the same for $I_{cb}$); $I_a$ the indicator of an *a*-vertex (the same for $I_b$). In (Gorbunov and Lyubetsky [21], Lemma 5), an important equality is proved:

**Claim 2**. *For any breakpoint graph G we have* $t(G) = B(G)+S(G)+D(G)–P(G)$ (compare Claims 1 and 2). □

**Lemma 1***a*. *For any breakpoint graph G the algorithm performs* $t(G)+t_1(G)$ *operations in total, where* $t_1(G)$ *is the number of interactions in Stage* 7 (*if it is invoked, otherwise* 0).

**Proof.** *a*) We have $t(G)= t_0(G)+t(G_1)= t_0(G)+(B+S+D–P)(G_1)$, where $t_0(G)$ is the *number* of operations in Stages 1–2 and $G_1$ is the corresponding graph after the final operation of Stages 1–2, but $P(G_1)=0$. Now consider finalizing $G_1$ by the algorithm at Stages 3–8. In these stages (except for the first operations of Stages 6.1 and 7), $(B+S+D)(G_1)$ reduces precisely by 1 till 0. Moreover, in the first operation of Stage 6.1, $B+S+D$ does not change, but in the second operation it is reduced by at once 2; and in the first operations of each point of Stage 7, $B+S+D$ does not change, but in the second operation it is reduced by only 1. Thus, the total number of steps in the algorithm is $t_0(G)+(B+S+D)(G_1)+t_1(G)$. Therefore, it is equal $t(G)+t_1(G)$. □

*Recall* that $w_1=w_a$ is the cost of removal of an *a*-vertex, and $w_2=w_b$ is the cost of removal of a *b*-vertex. Recall also that in all reasonings we may assume that $w=1$ and, accordingly, $w_a, w_b \geq 1$, denote $\varepsilon_a=w_a-1\geq0$, $\varepsilon_b=w_b-1\geq0$, $\varepsilon=\varepsilon_a+\varepsilon_b$; this follows from the premises of Theorem 1 above and Theorem 3 below. In each of the following equalities, the left-hand side is simply a *notation* for the right-hand side convenient for further presentation: $U_a=C_a+I_{ab}+I_{pa}$, $U_b=C_b+I_{ab}+I_{pb}$, $A_a=C_a–I_{ca}\cdot(1–I_{ab})$, $A_b=C_b–I_{cb}\cdot(1–I_{ab})$. Let us **define**

$T_1(G) = t(G)+\varepsilon_a\cdot U_a+\varepsilon_b\cdot U_b$, $T_2(G) = t(G)+\varepsilon_a\cdot U_a+\varepsilon_b\cdot I_b+A_b$, and $T_3(G) = t(G)+\varepsilon_a\cdot I_a+\varepsilon_b\cdot I_b+A_a+A_b$.

*Let*, for instance, $w_a \leq w_b$, and **consider three cases**:
  1) $1\leq w_a\leq w_b\leq 2$, 2) $1\leq w_a\leq 2<w_b$, and 3) $1\leq 2<w_a\leq w_b$.
In what follows, numbers 1–3 indicate these three cases.

**Lemma 1***b*. *The algorithm cost c(G) is one of the following numbers*: 1) $T_1$, $T_1+\varepsilon_a$, $T_1+\varepsilon_b$, $T_1+\varepsilon$; 2) $T_2$, $T_2+\varepsilon_a$, $T_2+1$, $T_2+\varepsilon_a+1$; 3) $T_3$, $T_3+1$, $T_3+2$.

**Proof.** The three cases are treated similarly; let us prove the first of them. In all stages preceding Stage 8, the Rem operation is not used; however, singular vertices can be removed as a result of joining them when performing one of the intermergings. At Stage 8, there remain no more than one $(a,b)$-cycle and no more than one $a$- or $b$-loop, as well as *some problem chains* of size 1 (in cases of $3a'$ and $3b'$) or 0 (otherwise). At Stage 8, removals are made by the Rem operation only from cycles of size 2 (with one or two singular vertices) and loops (with one singular vertex), their number being $C_a+C_b+2I_{ab}+I_{pa}+I_{pb}$ plus *at most two removals* from problem chains. The total number of $a$-removals is $U_a+n_a$, where $n_a$ is either 0 or 1, and the same for $b$-removals (ref. to the end of Stage 3). The total cost of all operations in the algorithm is $(1+\varepsilon_a)(U_a+n_a)+(1+\varepsilon_b)(U_b+n_b)+(B-U_a-U_b-n_a-n_b)+(S+D-P)$, where the third term is the number of all special operations at Stages 1–7. The fourth term is the number of nonspecial operations in the algorithm. Indeed, according to Lemma $1a$ this number is $B+S+D-P+t_1-B$, but $t_1=0$ for case 1. Using Claim 2, we continue the chain of equalities: $= t+\varepsilon_a(U_a+n_a)+ \varepsilon_b(U_b+n_b) = T_1+(\varepsilon_a n_a+\varepsilon_b n_b)$. □

**Theorem 3.** *Let* $\min(w_1,w_2) \geq w>0$. *In the all three cases, the minimum cost $c(G)$ of reduction of the graph G is, respectively, not less than $T_1$, $T_2$, and $T_3$.*

Theorem 3 and Lemma $1b$ immediately implies the following: the additive error of the algorithm is, respectively, not greater than $\varepsilon$, $\varepsilon_a+1$, and 2, which proves the main result of the paper, Theorem 1, since linear runtime of the algorithm easily follows from its description and Lemma $1a$. Indeed, the number of executed operations is $t(a+b)+t_1(a+b)$. The sum linearly depends on the size of the original structures $a$ and $b$. Note that to perform the next interaction, we need to have list of all components of a fixed type for a current graph $G$. Such list is obtained in constant time $C_1$ from the lists for the preceding $G_0$ and the transition operation from $G_0$ to $G$. The lists for the initial graph $a+b$ are composed by looking over all its components, which requires $C_2$ time linearly depending on the initial sizes of $a$ and $b$. Thus, the algorithm requires $C_1 \cdot (t(a+b)+t_1(a+b)) + C_2$ time.

We say that any breakpoint graph $G$ is of the $1^{st}$ *rank* if it has no $(a,b)$-components but contains $a$- and $b$-chains; and of the $2^{nd}$ *rank* if it has no $(a,b)$-components and contains either an $a$-chain or a $b$-chain but not both together (*named $a$- and $b$-graphs*). All the other graphs are of the $3^{rd}$ *rank*. Recall that an $a$-chain is an $a$-component that is a chain.

The inequality

$$c(o) \geq T_1(G) - T_1(o(G)) + d$$

is called the ***triangle inequality*** (for the first case); the additional term $d$ will be referred to as an *increment*. Specifically, $o$ transforms the $3^{rd}$ rank into the $2^{nd}$ (then the increments are $-\varepsilon_a$ and $-\varepsilon_b$ for $a$- and $b$-graphs) or vice versa (the increments are $+\varepsilon_a$ and $+\varepsilon_b$ for $a$- and $b$-graphs); or the $3^{rd}$ into the $1^{st}$ or vice versa (the increments are $-\varepsilon$ and $+\varepsilon$); or the $2^{nd}$ into the $1^{st}$ or vice versa (the increments are minus or plus $\varepsilon_b$ for an $a$-graph and minus or plus $\varepsilon_a$ for a $b$-graph); the *ordinary* triangle inequality means $d=0$.

**Proof of Theorem 3.** Consider the *first case* $1 \leq w_a \leq w_b \leq 2$. It suffices to check the triangle inequality and Claim 3 below, because $\varepsilon_a$ and $\varepsilon_b$ are nonnegative. □

**Claim 3.** *Assume that the <u>triangle inequalities</u> for the first case <u>hold</u>. The minimum cost $C(G)$ of reduction of G for a graph of the $1^{st}$ rank is $\geq T_1(G)+\varepsilon$; of the $2^{nd}$ rank, either $\geq T_1(G)+\varepsilon_a$ or $\geq T_1(G)+\varepsilon_b$, respectively; and of the $3^{rd}$ rank, $\geq T_1(G)$.*

**Proof.** Consider the first case. Using induction on the minimum cost $C(G)$ of reduction of $G$, we show that these inequalities is implied by the corresponding triangle inequality. For instance, if $o$ takes $G$ of the 1st rank to an $a$-graph of the 2nd rank, by the induction hypothesis we obtain $C(o(G)) \geq T_1(o(G)) + \varepsilon_a \geq T_1(G) - c(o) + \varepsilon_a + \varepsilon_b$ and then $C(G) = C(o(G)) + c(o) \geq T_1(G) - c(o) + \varepsilon + c(o) = T_1(G) + \varepsilon$; here, $o$ is the first operation in the shortest transformation for $G$.

When passing between graphs of the same rank, the desired inequality follows from the ordinary triangle inequality. Transitions between an $a$- and a $b$-graph of the 2nd rank are impossible. The other transitions are considered similarly. □

Now comes a significant statement.

**Lemma 2.** *Consider the first case. The triangle inequalities are valid.*
**Proof.** Consider all operations Rem, OM, Cut, SM and DM. We start with the operations **Rem** for an $a$-vertex; removal of a $b$-vertex is treated similarly. If an $a$-vertex is removed from a cycle or loop, the rank of the graph does not increase, $U_b$ does not decrease (if $I_{ab}$ decreases by 1, then $C_b$ increases by 1), and $U_a$ can decrease by at most 1 (since exactly one of the numbers $C_a$, $I_{ab}$, and $I_{pa}$ can decrease). The triangle inequality follows from the $t$-property.

Let a vertex be removed from a chain. If $I_{ab}$ does not decrease, the rank either does not change (then the result is implied by the $t$-property and $\varepsilon_a \geq 0$) or changes from the 1st to the $b$-second or from the $a$-second to the 3rd (then adding $\varepsilon_a$ to 1, i.e., to the maximum possible reduction of $t$, gives exactly $w_a$). If $I_{ab}$ decreases, then $T_1$ decreases by at most $1+\varepsilon$. Then the rank changes either from the 3rd to the 1st (then subtract $\varepsilon$ from $1+\varepsilon$) or from the 3rd to the $b$-second (then subtract $\varepsilon_b$ from $1+\varepsilon$ and obtain $1+\varepsilon_a$). Before continuing, let us prove the following.

**Claim 4.** *If after operations OM, SM, or DM (denoted by $o$) a loop is formed, then $t(G) \leq t(o(G))$.*
**Proof.** Assume that $t(o(G)) = t-1$, where $t = t(G)$, and use Claim 2. After the Rem operation over the loop in $o(G)$, we have $t(G_1) = t-2$, since $B$ decreases by 1 and the other characteristics do not change. The same graph $G_1$ is obtained by the Rem operation over a singular vertex from which the loop was formed. This contradicts the $t$-property. Note that after the Rem operation, changes in characteristics of $G$ are not so clear. □

Let us continue the proof of Lemma 2. Consider an **OM** operation $o$. Then $U_a$ and $U_b$ do not decrease. If $o$ does not increase the rank, the result follows from the $t$-property. If, upon circularizing an $a$-chain (similarly, a $b$-chain), the rank changes from the 1st to the 2nd or from the 2nd to the 3rd, then $C_a$ increases by 1, and the maximum possible decrease of $T_1$ is $1-\varepsilon_a$. Adding $\varepsilon_a$ to this quantity precisely gives the cost of the OM. If the same increase of the rank occurs because of forming a loop, then the triangle inequality is implied by Lemma 2 and $\varepsilon_a \leq 1$. If the rank changes from the 1st to the 3rd ($a$- and $b$-chains are merged), then $I_{ab}$ increases by 1, and the maximum possible decrease of $T_1$ is $1-\varepsilon$. Adding $\varepsilon$ to this quantity gives 1.

Consider a **Cut** operation. It does not increase the rank, and taking into account the $t$-property it suffices to consider cases where $U_a$ or $U_b$ decreases. Let an $a$-loop be cut. Since $I_{pa}$ decreases, the rank changes from the 3rd to the $a$-second or from the $b$-second to the 1st. The decrease of $T_1$ is at most $1+\varepsilon_a$. Subtracting $\varepsilon_a$ gives 1. Let an $a$-cycle be cut. Then $C_a$ decreases but $t$ does not; see (Gorbunov and Lyubetsky [23], item 3 in the proof of Theorem 2). Then $T_1$ decreases by at most $\varepsilon_a \leq 1$. Let an $(a,b)$-chain be cut into an $a$-chain and a $b$-chain. $I_{ab}$ decreases, and the rank changes from the 3rd to the 1st. The decrease of $T_1$ is at most $1+\varepsilon$. Subtracting $\varepsilon$ gives 1. The same for a $b$-loop and an $a$-cycle.

Consider a **SM** operation. Let the rank increase. If it increases after transforming an $a$-chain or a $b$-chain into a cycle or loop, we repeat the arguments for an OM. If the rank increases because of forming an $(a,b)$-component from an $a$-component and a $b$-component, three cases are possible. In the first case, an $a$-chain and a $b$-chain are merged; then we repeat the arguments for an OM. In the second case, assume that an $a$-cycle and a $b$-chain are merged; then $I_{ab}$ increases, $C_a$ decreases, the rank changes from the 1st (or $b$-second) to 3rd. Then the triangle inequality follows from the fact that $t$ does not decrease under a SM with cutting an $a$-cycle or a $b$-cycle; see (Gorbunov and Lyubetsky [23], item 4 in the same proof). In the third case, assume that an $a$-loop and a $b$-chain are merged; then $I_{ab}$ increases. If $I_{pa}$ decreases, the rank changes from the $b$-second to the 3$^{rd}$. The possible decrease of $T_1$ is $1-\varepsilon_b$; adding $\varepsilon_b$ to it gives 1. If $I_{pa}$ does not decrease, the rank changes from the 1$^{st}$ to the 3$^{rd}$. The possible decrease of $T_1$ is $1-\varepsilon$; adding $\varepsilon$ gives the same effect.

Assume that the rank does not increase. It suffices to consider the cases where $U_a$ or $U_b$ decreases. If $C_a$ or $C_b$ decreases, then $t$ does not decrease (the same item 4 of Theorem 2 in Gorbunov and Lyubetsky [23]). Then $T_1$ decreases by at most $\varepsilon_a \leq 1$. Assume that $I_{ab}$ decreases, i.e., an $(a,b)$-chain is cut into an $a$-chain and a $b$-chain. If neither of them is circularized into a cycle or loop, repeat the arguments for a cut. Otherwise: let, for instance, an $a$-chain be circularized. Then the rank changes from the 3$^{rd}$ to the 1$^{st}$ or $b$-second. In the case of a cycle, $C_a$ increases, then $T_1$ decreases by $1+\varepsilon_b$ and the increment is $\varepsilon_b$, and in the case of a loop the result follows from Claim 4. Assume that $I_{pa}$ or $I_{pb}$ (for definiteness, $I_{pa}$) decreases, i.e., an $a$-loop is cut off and inserted into a chain. Then $I_{ab}$ does not increase, since otherwise the rank would change from the $b$-second to the 3$^{rd}$. Repeat the arguments for a cut.

Consider a **DM** operation. For it, the $t$-property is refined in Claim 5. A maximal connected chain of conventional edges is called a *segment*. In the proof of Claim 5 we will use the following observation. Cutting off all conventional edges from a segment of length $l$ located between two singular edges in a cycle of size strictly greater than 2 requires $l/2$ nonspecial operations if $l$ is even and $(l-1)/2$ nonspecial operations otherwise.

**Claim 5.** *The result of applying a DM operation that joins an $a$- or $b$-cycle with $b$-loop or $a$-loop, respectively, does not reduce $t(G)$. If an $a$-cycle is joined with a $b$-cycle, then $t$ increases by 1.*

**Proof.** Let an $a$-cycle be joined with a $b$-loop; similarly for a $b$-cycle. Between $a$-vertices in an $a$-cycle there is an odd number of conventional edges. The first statement follows from the observation that $B$, $S$, $D$, and $P$ are unchanged. Indeed, here the DM consists in removing a loop and replacing a conventional $b$-edge in an $a$-cycle wih two singular edges. As a result, a segment of an odd length is divided into two segments of even lengths (to the empty segment, we assign length 0) of the total length less by 1 than the length of the original segment. The second statement follows from the observation that $B$, $D$, and $P$ are unchanged and two odd segments of conventional edges turn into two segments such that each of them is located between an $a$-vertex and a $b$-vertex, and therefore has an even length. □

Using Claim 5, we check Lemma 2 for DM analogously to SM. This completes the proof of Lemma 2 <u>in the first case</u>. □

Now consider the *second case* $1 \leq w_a \leq 2 < w_b$. For that, we change the definition of the rank of a graph. We say that a breakpoint graph $G$ is of the 1$^{st}$ *rank* if it has no $(a,b)$-components but contains $a$-chains, $b$-chains, and $b$-cycles; and of the 2$^{nd}$ *rank* if it has no $(a,b)$-components and contains either an $a$-chain or a $b$-chain and a $b$-cycle but not both together ($a$-graphs or $b$-graphs, respectively); asymmetry in this definition is followed from asymmetry of $w_a$ and $w_b$. All the other graphs are of the 3$^{rd}$ *rank*. Claim 3 takes the following form: the minimum cost $C$ of reduction of $G$ for a graph of the 1$^{st}$ *rank* is $\geq T_2+\varepsilon_a+1$; of the 2$^{d}$ *rank*, either $\geq T_2+\varepsilon_a$ or $\geq T_2+1$, respectively; and of the 3$^{d}$ *rank*, $\geq T_2$. □

Now consider the *third case* $2 < w_a \leq w_b$. For that, we again change the definition of the rank. We say that a breakpoint graph is of the 1$^{st}$ *rank* if it has no $(a,b)$-components but contains $a$-chains, $a$-cycles, $b$-chains, and $b$-cycles, and of the 2$^{d}$ *rank* if it has no $(a,b)$-components and contains either an $a$-chain and an $a$-cycle or a $b$-chain and a $b$-cycle but not both together ($a$-graphs and $b$-graphs). All the other graphs are of the 3$^{rd}$ *rank*. Claim 3 takes the following form: the minimum cost $C$ of reduction of $G$ for a graph of the *first rank* is $\geq T_3+2$; of the *second* *rank*, $\geq T_3+1$; and of the *third rank*, $\geq T_3$. □

After these changes, the arguments in the second and third cases are similar to those in the first case. This completes the proof of Lemma 2. □

**Corollary of Theorem 3.** *If* $1=w_a \leq w_b \leq 2$, *then the additive error of the algorithm is at most* $w_b-1$. *If all components in initial structures are circular, then the algorithm is exact.*

**Acknowledgments** The research was supported by the Russian Foundation for Basic Research under research project No. 18-29-13037.

**References**
1. Bioinformatics and Phylogenetics. Seminal Contributions of Bernard Moret. Eds. Tandy Warnow. Comput. Biol. Series. Springer Nature Switzerland AG, 2019.
2. Models and Algorithms for Genome Evolution. Eds. C. Chauve, N. El-Mabrouk, E. Tannie. Comput. Biol. Series. London: Springer, 2013
3. Pavel A. Pevzner. Computational Molecular Biology: An Algorithmic Approach. The MIT Press, Cambridge, Massachusetts, London, England, 2000.
4. D. Sankoff, G. Leduc, N. Antoine, B. Paquin, B. F. Lang, R. Cedergren (1992). Gene order comparisons for phylogenetic inference: evolution of mitochondrial genome. Proc Natl Acad Sci **89**:6575–6579.
5. S. Hannenhalli, P. Pevzner (1999). Transforming cabbage into turnip: polynomial algorithm for sorting signed permutations by reversals. J of the ACM **46**(1):1–27.
6. S. Hannenhalli, P. Pevzner (1995). Transforming man into mice (polynomial algorithm for genomic distance problem). Proc. of the 36th Annual Symposium on Foundations of Computer Science, p. 581.
7. M. A. Alekseyev, P. A. Pevzner (2008). Multi-Break Rearrangements and Chromosomal Evolution. Theor Computer Science **395**(2–3): 193–202.
8. M. D. V. Braga, E. Willing, J. Stoye (2011). Double cut and join with insertions and deletions. J. Comput. Biol. **18**(9):1167–1184.
9. P. E. C. Compeau (2013). DCJ-indel sorting revisited. Algorithms for Mol Biol **8**:6.
10. P. H. da Silva, R. Machado, S. Dantas and M. D. V. Braga (2013). DCJ-indel and DCJ-substitution distances with distinct operation costs. Algorithms for Mol Biol **8**:21.
11. P. E. G. Compeau (2014). A Generalized Cost Model for DCJ-Indel Sorting. LNCS, vol. 8701, pp 38–51.
12. P. H. da Silva, R. Machado, S. Dantas and M. D. V. Braga (2017). Genomic Distance with High Indel Costs. J IEEE/ACM Trans Comput Biol Bioinform **14**(3): 1–6.
13. H. T. Jacobs, D. J. Elliot, V. B. Math and A. Farquharson (1988). Nucleotide sequence and gene organization of sea urchin mitochondrial DNA. J Mol Biol **202**(2):185–217.
14. R. H. Bors, J. A. Sullivan (2005). Interspecific Hybridization of *Fragaria vesca* subspecies with *F. nilgerrensis*, *F. nubicola*, *F. pentaphylla*, and *F. viridis*. Journal of the American Society for Horticultural Science **130**(3):418–423.
15. N. H. Putnam, T. Butts, D. E. K. Ferrier, R. F. Furlong, U. Hellsten, T. Kawashima, M. Robinson-Rechavi, E. Shoguchi, A. Terry, Jr-Kai Yu, È. Benito-Gutiérrez, I. Dubchak, J. Garcia-Fernàndez, J. J. Gibson-Brown, I. V. Grigoriev, A. C. Horton, P. J. de Jong, J. Jurka, V. V. Kapitonov, Y. Kohara, Y. Kuroki, E. Lindquist, S. Lucas, K. Osoegawa, L. A. Pennacchio, A. A. Salamov, Y. Satou, T. Sauka-Spengler, J. Schmutz, T. Shin-I, A. Toyoda, M. Bronner-Fraser, A. Fujiyama, L. Z. Holland, P. W. H. Holland, N. Satoh & D. S. Rokhsar (2008). The amphioxus genome and the evolution of the chordate karyotype. Nature **453**:1064–1072.
16. M. Shao, Y. Lin, B. Moret (2014). An exact algorithm to compute the DCJ distance for genomes with duplicate genes. In: Proc. of RECOMB 2014, LNBI, vol. 8394, pp. 280–292, Heidelberg: Springer Verlag.

17. F. V. Martinez, P. Feijão, M. D. V. Braga and J. Stoye (2015). On the family-free DCJ distance and similarity. Algorithms for Mol Biol **10**:13.
18. V. A. Lyubetsky, R. A. Gershgorin, A. V. Seliverstov, K. Yu. Gorbunov (2016). Algorithms for Reconstruction of Chromosomal Structures. BMC Bioinformatics **17**:40.1–40.23.
19. V. A. Lyubetsky, R. A. Gershgorin, K. Yu. Gorbunov (2017). Chromosome Structures: Reduction of Certain Problems with Unequal Gene Content and Gene Paralogs to Integer Linear Programming. BMC Bioinformatics **18**:537.1–537.18.
20. A. Bergeron, J. Mixtacki, J. Stoye (2006). A unifying view of genome rearrangements. Algorithms in Bioinformatics, LNCS, vol. 4175, pp. 163–173.
21. K. Yu. Gorbunov, V. A. Lyubetsky (2017). Linear algorithm of the minimal reconstruction of structures. Probl of Inform Transmission **53**(1):55–72.
22. S. Yancopoulos and R. Friedberg (2009). DCJ Path Formulation for Genome Transformations which Include Insertions, Deletions, and Duplications. J. of Comput. Biology **16**(10): 1311–1338.
23. K. Yu. Gorbunov, V. A. Lyubetsky (2017). A linear algorithm for the shortest transformation of graphs with different operation costs. J of Commun Technology and Electronics **62**(6):653–662.